\documentclass[a4paper,fontsize=14pt]{scrartcl}
\usepackage{amsmath,amssymb,amsthm}
\usepackage{microtype}
\usepackage{hyperref}

\newtheorem{theorem}{Theorem}
\newtheorem{lemma}[theorem]{Lemma}
\newtheorem{corollary}[theorem]{Corollary}
\newtheorem{proposition}[theorem]{Proposition}
\theoremstyle{definition}
\newtheorem{definition}[theorem]{Definition}
\theoremstyle{remark}
\newtheorem{remark}{Remark}
\newtheorem{example}{Example}

\newcommand\lcm{\mathrm{lcm}}
\newcommand\diam{\mathrm{d}}
\newcommand\FF{\mathbb{F}}
\newcommand\FFp{\FF_p}
\newcommand\intrv[1]{{[#1]}}

\title{Generalizing the Bierbrauer--Friedman bound for % mixed- and pure-level
orthogonal arrays%
\thanks{This is a post-peer-review, pre-copyedit version of the article published in Designs, Codes and Cryptography 93(11), 2025, 4937--4950, \url{https://doi.org/10.1007/s10623-025-01711-y}}%
\thanks{The work of D.K. was supported by the Natural Science Foundation of Hebei Province (project No. A2023205045)
and within the framework of
the state contract of the Sobolev Institute of Mathematics (FWNF2022-0017).\\
The work of F.\"O. was supported by T\"UB\.{I}TAK under Grant 223N065.
}%
}
\author{Denis S. Krotov%
\thanks{School of Mathematical Sciences, Hebei Key Laboratory of Computational Mathematics and Applications, Hebei Normal University, Shijiazhuang 050024, P. R. China}
\thanks{Sobolev Institute of Mathematics, Novosibirsk 630090, Russia%; \texttt{https://orcid.org/0000-0002-8516-755X}
}%
, \
Ferruh \"Ozbudak%
\thanks{Faculty of Engineering and Natural Science, Sabanc{\i} University, Istanbul, Turkey; 
email: \texttt{ferruh.ozbudak@sabanciuniv.edu}}%
, \
Vladimir N. Potapov%
\thanks{Faculty of Engineering and Natural Science, Sabanc{\i} University, Istanbul, Turkey%;
}%
}
\date{}

\begin{document}

\maketitle

\begin{abstract}
We characterize mixed-level orthogonal arrays in terms of algebraic designs in a special multigraph. We prove a mixed-level analog of the Bierbrauer--Friedman (BF) bound for pure-level orthogonal arrays and show that arrays attaining it are radius-$1$ completely regular codes (equivalently, intriguing sets, equitable $2$-partitions, perfect $2$-colorings) in the corresponding multigraph. For the case when the numbers of levels are powers of the same prime number, we characterize, in terms of multispreads, additive mixed-level orthogonal arrays attaining the BF bound. For pure-level orthogonal arrays, we consider versions of the BF bound obtained by replacing the Hamming graph by its polynomial generalization and show that in some cases this gives a new bound.

\end{abstract}

\section{Introduction}
Orthogonal arrays are combinatorial structures important both for practical
applications like design of experiments or software testing and for theoretical purposes,
because of many relations with coding theory, cryptography, design theory,
etc., see e.g. \cite{HSS:OA}.
Among many other interesting relations,
pure-level, or symmetric, orthogonal arrays
are known as a special case
of algebraic designs,
which makes them a part of a general theory
that includes also other widely known classes of
combinatorial objects, such as combinatorial $t$-$(v,k,\lambda)$ designs.
One of the main results of this correspondence
is establishing a similar relation for mixed-level (asymmetric) orthogonal arrays,
attracting more attention in recent years, see e.g.
recent works
\cite{ChenNiu:2023}, \cite{LinPangWang:2024}, \cite{NiuChenGaoPang:2025} and references therein. Then,
we use the correspondence obtained
to generalize results known for pure-level
orthogonal arrays,
namely, the Bierbrauer--Friedman bound
and constructions of arrays attaining it.
Additionally, we consider some generalized
versions of this bound for pure-level arrays
that give a nontrivial inequality
when the original bound is not applicable.

The Bierbrauer--Friedman bound
for pure-level orthogonal arrays with parameters
OA$(N,n,q,t)$ says that
\begin{equation}\label{eq:BFpu}
N \ge q^n \Big(1-\big(1-\frac1q\big)\frac{n}{t+1}\Big),
\end{equation}
see~\cite{Friedman:92} for
the case $q=2$ and~\cite{Bierbrauer:95} for general~$q$.
It is easy to see that the bound
is nonnegative if and only if
$t+1 > \frac{q-1}{q}n$; so, it is effective for high values of~$t$
(in contrast, Rao's bound~\cite{Rao:47} is effective for relatively small~$t$).
The bound is tight,
and there are orthogonal arrays
constructed as linear \cite[Section\,4.3]{HSS:OA}
or additive \cite[Section\,4.2]{BKMTV} codes that attain it.
Binary ($q=2$)
orthogonal arrays attaining
bound~\eqref{eq:BFpu}
whose size is not a power
of~$2$ can be constructed
as completely regular codes
(see Definition~\ref{def:CR} below)
by the Fon-Der-Flaass construction~\cite{FDF:PerfCol};
the first example is OA$(1536, 13, 2, 7)$.
Ternary ($q=3$) arrays with
the similar property were recently
discovered, again in terms of completely regular codes, in~\cite{GorKro:CR};
the first example is OA$(5\cdot3^8, 11, 3, 8)$.
%% TODO: add reference to the table
For $q\ge 4$, the problem of existence of orthogonal arrays attaining bound~\eqref{eq:BFpu} whose size
is not a power of~$q$ or of its prime divisor
remains open. Similar questions can be considered
for mixed-level orthogonal arrays.

In this correspondence,
we prove (Section~\ref{s:mixed}) that \eqref{eq:BFpu}
holds for mixed-level
orthogonal arrays
OA$(N,q_1\cdot q_2\cdot\ldots\cdot q_n,t)$
if we replace $q^n$
by the product of all $q_i$s
and $\frac1q$
by the average value
of~$\frac{1}{q_i}$.
The arrays (we treat an array as a multiset of rows of length~$n$) attaining this bound
are necessarily simple sets
(without repeated elements),
independent sets
(without pairs of elements
at Hamming distance~$1$),
and intriguing sets
(completely regular codes
with covering radius~$1$,
see the definition below).
The new bound is tighter than
the previously known generalization~\cite{Diestelkamp:2004}
\def\Max{\mathrm{M}}
\def\Min{\mathrm{m}}
\begin{align}
\label{eq:BFD}
N \ge q_{\Min}^n \Big(1 -
\frac{n\tilde q-n}{n\tilde q + (t+1-n)q_{\Max}}\Big)
% \\[1ex] \nonumber
\qquad \mbox{if}\ \  {n\tilde q + (t+1-n)q_{\Max}}>0,
\end{align}
where $q_{\Min}$, $\tilde q$,
and $q_{\Max}$ are
respectively the minimum,
the average, and the maximum value of~$q_i$, $i=1,\ldots,n$.
For example, for
OA$(N,2^1 4^4,3)$
% OA$(N,2{\cdot} 4{\cdot} 4{\cdot} 4{\cdot} 4,3)$
the new bound is tight:
$N \ge 64 = 2^9\cdot\big(1-(1-\frac{3}{10})\frac{5}{3+1}\big)$
(see the construction in Example~\ref{ex:1}),
while
\eqref{eq:BFD} gives
$N \ge \frac{16}{7}=2^5\cdot\big(1-\frac{(3.6-1)\cdot 5}{4\cdot(3+1)-(4-3.6)\cdot 5}\big)$.

Further (Section~\ref{s:additive}), in the case when all $q_i$
are powers of the same prime~$p$,
we prove that additive (linear over GF($p$))
mixed-level orthogonal arrays attaining the Bierbrauer--Friedman bound
are equivalent to special partitions of
a vector space into subspaces, called multispreads~\cite{KroMog:multispread}.

Finally (Section~\ref{s:PG-BFB}), we discuss variations of
bound~\eqref{eq:BFpu}
that can give positive values
when the original bound~\eqref{eq:BFpu}
is negative. As an example,
for OA$(N,n,2,\frac n2-1)$, $n$ even,
we obtain the bound
$ N > {0.409}\cdot {n}^{-1} \cdot 2^n$.
To compare, for these parameters, Rao's bound
gives the size of the Hamming ball
of radius approximately
$n/4$,
which is
$2^{h(\frac14)n(1+o(1))}$,
$h(\frac14)\simeq 0.8$.

Our main results are
Theorem~\ref{th:OAtD}
(mixed-level orthogonal arrays are algebraic
designs),
Theorem~\ref{th:gBF}
(Bierbrauer--Friedman bound for mixed-level orthogonal arrays and a relation with completely regular codes),
Theorem~\ref{th:additive}
(the characterization of
additive mixed-level orthogonal arrays attaining the Bier\-brauer--Friedman bound),
and Lemma~\ref{l:gBF1}
(a polynomial generalization of
the Bier\-brauer--Friedman bound
for mixed-level orthogonal arrays).

\section{Definitions and notations}\label{s:def}

By a graph, we will mean a multigraph, with multiple edges and loops allowed.
A graph without loops and edge multiplicities more than~$1$ is called simple.
For a graph~$G$ and a positive integer scalar~$s$, $sG$ denotes the graph
on the same vertex set with all edge multiplicities multiplied by~$s$.
For two graphs~$G'=(V',E')$
and~$G''=(V'',E'')$, $G'\square G''$
denotes their Cartesian product,
the graph with the vertex set
$V'\times V''$ and the edge (multi)set
$\{ \{(v_1,v),(v_2,v)\}: \{v_1,v_2\}\in E',\, v\in V'' \}
\cup
\{ \{(v,v_1),(v,v_2)\}: \{v_1,v_2\}\in E'',\, v\in V' \}$.

For a positive integer $s$, the set
$\{0,\ldots,s-1\}$ will be denoted by~$\intrv{s}$.
\begin{definition}[adjacency matrix of a graph, eigenfunctions and eigenspaces]
For a graph $G=(V,E)$, the \emph{adjacency matrix} $A$ is the symmetric nonnegative integer  $|V|\times |V|$ matrix whose rows and columns are indexed by~$V$ and the $(x,y)$th element $A_{x,y}$, $x,y\in V$,
equals the multiplicity of $\{x,y\}$ in the edge multiset~$E$.
The eigenvectors of~$A$, treated as functions from~$V$ to~$\mathbb{C}$, are called \emph{eigenfunctions} of~$G$, or \emph{$\theta$-eigenfunctions}, where $\theta$ is the corresponding eigenvalue.
An \emph{eigenspace} of~$G$
is the subspace of the vector space
$\mathbb{C}^{V}$
consisting of the constantly zero function and all  $\theta$-eigenfunctions
for some eigenvalue~$\theta$.
Since the adjacency matrix of a graph is symmetric,
the eigenspaces are pairwise orthogonal.
\end{definition}

\begin{definition}[completely regular code with covering radius~$1$, $\{b;c\}$-CR code] \label{def:CR}
A set $C$ of vertices
of a regular graph $G=(V,E)$
is called a \emph{completely regular code} with covering radius~$1$ and intersection array $\{b;c\}$, where $b,c>0$,
or a \emph{$\{b;c\}$-CR code}, or  simply a \emph{CR-$1$ code}
if for every vertex in~$V$ (in~$V\backslash C$)
the number of edges that connect it
with~$V\backslash C$
(respectively, with~$C$),
equals~$b$ (respectively,~$c$).
\end{definition}

\begin{remark} In literature, CR-$1$ codes are studied under different names. They are also called
\emph{intriguing sets}; the corresponding
partition of the vertex set into $C$ and $V\backslash C$ is known as an \emph{equitable $2$-partition} or a \emph{$2$-partition design};
the corresponding $2$-coloring of the vertex set is a \emph{perfect $2$-coloring}.
\end{remark}

\begin{definition}[$H(q_1\cdot q_2\cdot\ldots\cdot q_n)$, Hamming graph]
For integers $q_1,q_2,...,q_n\ge 2$,
the graph
$H(q_1\cdot q_2\cdot\ldots\cdot q_n)$ is the graph on the set
of $n$-tuples from
$V=\intrv{q_1}\times \ldots \times \intrv{q_n}$
with
two $n$-tuples forming an edge
of multiplicity~$\mu$ if and only if
they differ in only the $i$th position for some $i\in\{1,\ldots,n\}$ and $\mu = Q/q_i$,
where $Q=\lcm(q_1,...,q_n)$ (least common multiple).
If $q_1=q_2=\ldots=q_n=q$,
then $H(q_1\cdot q_2\cdot\ldots\cdot q_n)$ is a simple graph known as a \emph{Hamming graph} and denoted $H(n,q)$.

In particular, $H(1,q)$ is the complete graph of order~$q$ and
$$ H(q_1 \cdot\ldots\cdot q_n) =
\frac{Q}{q_1}H(1,q_1)\,\square\, \ldots \,\square\, \frac{Q}{q_n}H(1,q_n).  $$
\end{definition}

\begin{definition}[algebraic $t$-design]
For a regular graph $G=(V,E)$ with  eigenvalues $\theta_0>\theta_1>\ldots>\theta_d$
and the corresponding eigenspaces
$S_0$, $S_1$, \ldots, $S_d$,
a multiset $C$ of its vertices is called
an \emph{algebraic $t$-design} (with respect to the natural descending ordering  $\theta_0,\theta_1,\ldots,\theta_d$ of the graph eigenvalues) if in the decomposition
$$
f_C = \varphi_0+\varphi_1+\ldots+\varphi_d,
\qquad \varphi_i \in S_i,
$$
of the multiplicity function~$f_C$ (in the case of a simple set,
the characteristic function, indicator) of~$C$
we have $\varphi_1=\ldots=\varphi_t\equiv 0$.
\end{definition}

\begin{definition}[orthogonal arrays, OA$(N,q_1{\cdot} q_2{\cdot} \ldots{\cdot} q_n , t)$]
A nonempty multiset~$C$ of $n$-tuples
from
$V=\intrv{q_1}\times \ldots \times \intrv{q_n}$
is called an \emph{orthogonal array}
of strength~$t$,
OA$(|C|,q_1\cdot q_2 \cdot \ldots\cdot  q_n,t)$, if for any
distinct $i_1$, \ldots, $i_t$ from
$\{1,\ldots,n\}$ and any
$a_1\in q_{i_1}$, \ldots,
$a_t\in q_{i_t}$,
the number of
$(x_1,...,x_n)\in C$ such that
$x_{i_j}=a_j$, $j=1,...,t$,
equals $\frac{|C|}{q_{i_1}q_{i_2}...q_{i_t}}$ (i.e., independent on the choice of $a_j$, $j=1,...,t$).
If $q_1= q_2= \ldots = q_n$, then such orthogonal arrays are
called \emph{pure-level} (in some literature, \emph{symmetric}) and also denoted
OA$(N,n,q_1,t)$;
otherwise, they are called \emph{mixed-level}, or just \emph{mixed} (in some literature, \emph{asymmetric}).
For brevity, in the notation
$q_1\cdot q_2 \cdot \ldots \cdot q_n$,
equal values of $q_i$ can be grouped
using degrees, e.g.,
OA$(N,2{\cdot} 2 {\cdot} 5,t)$
is the same as OA$(N,2^2 5^1,t)$,
but not the same as
OA$(N,4 {\cdot} 5,t)= \text{OA}(N,4^1 5^1,t)$.
\end{definition}
%\begin{remark}\label{r:q} The definitions above allows to use different notations for the same parameters. For example, $OA(N,2^3 3^2,t)$ is a compact form of $OA(N,2^1 2^1  2^1 3^1 3^1,t)$. The last form is convenient to use in some proofs. Similarly, the graphs $H(2^3 3^2)$ and $H(2^1 2^1 2^1 3^1 3^1 )$ are the same. \end{remark}

\section{Orthogonal arrays and algebraic designs}\label{s:OA-AD}\label{s:mixed}

We first describe the eigenspaces of the graph
$H(q_1{\cdot}q_2{\cdot}...{\cdot}q_n)$.

\begin{lemma}
The following functions
form an orthogonal basis
from eigenfunctions of
$H(q_1\cdot q_2 \cdot \ldots \cdot q_n)$:
\begin{equation}\label{eq:bss}
\chi_{(b_1,b_2,...,b_n)}(x_1,x_2,...,x_n)
=
\xi_1 ^{b_1 x_1}
\xi_2 ^{b_2 x_2}
\ldots
\xi_s ^{b_n x_n},\qquad b_i\in [q_i], %\{0,\ldots,q_i-1\},
\end{equation}
where $\xi_i$ is the degree-$q_i$
primitive root of~$1$. Moreover,
$\chi_{(b_1,b_2,...,b_n)}$
is an eigenfunction corresponding to the
eigenvalue $\theta_w = k - w Q$,
where
$$
k = \sum_{i=1}^n
\frac{Q}{q_i}
(q_i-1)
$$
is the degree of $H(q_1{\cdot}q_2{\cdot}...{\cdot}q_n)$, $Q=\lcm(q_1,...,q_n)$,
and $w$ is the number of nonzeros
 among $b_1$, $b_2$, \ldots, $b_n$.
\end{lemma}
\begin{proof}
The complete graph $K_{q_i}$ has eigenfunctions
$\xi_{(b_i)}(x_i) = \xi_i^{b_ix_i}$
corresponding to the eigenvalue $q_i-1$
if $b_i=0$ and $-1$ otherwise.
After multiplying (the multiplicity of all edges of) $K_{q_i}$ by~$\frac{Q}{q_i}$, we get
the same eigenfunctions with eigenvalues
$\frac{Q}{q_i}(q_i-1)$ and $-\frac{Q}{q_i}(q_i-1)$.
We note that the difference~$Q$
between these two eigenvalues
does not depend on~$i$.

The rest is straightforward
from the following easy property
of the Cartesian product of graphs,
see e.g.
\cite[Section\,1.4.6]{BroHae:spectra}:
{if $\chi'(x')$, $x'\in V'$,
and $\chi''(x'')$, $x''\in V''$,
are $\theta'$- and $\theta''$-eigenfunctions
of graphs $\Gamma'=(V',E')$
and $\Gamma''=(V'',E'')$, respectively,
then $\chi(x',x'')=\chi'(x')\chi''(x'')$
is a $(\theta'+\theta'')$-eigenfunction
of~$\Gamma'\square\Gamma''$.}
%
%The claim is trivial if $n=1$, in which case we have the complete graph $K_{q_1}$ of order~$q_1$. In general, $H(q_1{\cdot}q_2{\cdot}...{\cdot}q_n) = \frac{Q}{q_1}K_{q_1} \square \ldots \square \frac{Q}{q_n}K_{q_n}$, and the claim follows from properties of eigenspaces of the Cartesian product of graphs (see e.g. \cite[Section\,1.4.6]{BroHae:spectra}) and simple arithmetic calculations.
\end{proof}

The following fact generalizes Delsarte's
characterization \cite[Theorem~4.4]{Delsarte:1973} of pure-level
orthogonal arrays as algebraic designs.
\begin{theorem}\label{th:OAtD}
A multiset $C$ of words
from $V=\intrv{q_1} \times \ldots \times \intrv{q_n}$
is an orthogonal array OA$(|C|, q_1{\cdot}q_2{\cdot}...{\cdot}q_n, t)$
if and only if $C$ is an algebraic
$t$-design in
$H(q_1{\cdot}q_2{\cdot}...{\cdot}q_n)$.
\end{theorem}
\begin{proof}
\emph{If.}
Assume $C$ is an algebraic
$t$-design and $f_C$ is its multiplicity function.
We need to show that
the sum of $f_C$
over $V_{i_1 ,\ldots, i_t}^{a_1 ,\ldots, a_t}$
does not depend on the choice of
$a_1$, \ldots, $a_t$,
where
$1\le i_1 < \ldots < i_t \le n$,
$a_j\in\intrv{q_{i_j}}$,
and $V_{i_1 ,\ldots, i_t}^{a_1 ,\ldots, a_t}$ denotes the set
of
all $n$-tuples
$(x_1,...,x_n)$ from~$V$
such that
$x_{i_1}=a_1$,
\ldots,
$x_{i_t}=a_t$.

We have
$$f_C = \varphi_0+\varphi_{t+1}+\varphi_{t+2}+\ldots+\varphi_{n},$$ where
$\varphi_i$ is the zero constant or a $\theta_i$-eigenfunction
 of $H(q_1{\cdot}q_2{\cdot}...{\cdot}q_n)$.

 For $i\ge t+1$, each such eigenfunction
 is a linear combination of
 basis eigenfunctions $\chi_{(b_1, \ldots, b_n)}$ from~\eqref{eq:bss}, where
 the number of nonzero elements among $b_1$, \ldots, $b_n$ is~$i$.
 We claim that
 \begin{itemize}
 \item[(*)] \emph{for more than $t$ nonzeros in $(b_1, \ldots, b_n)$, the sum of $\chi_{(b_1, \ldots, b_n)}$ over
 $V_{i_1 ,\ldots, i_t}^{a_1 ,\ldots, a_t}$ equals~$0$}. Indeed,
 denoting by $l_1$, \ldots, $l_{s}$ the indices from $\{1,...,n\}\backslash\{i_1,...,i_t\}$, $s=n-t$,
 we have
 \begin{multline}\label{eq:sums}
  \sum_{(x_1,...,x_n)\in V_{i_1 ,\ldots, i_t}^{a_1 ,\ldots, a_t}}
  \chi_{(b_1, \ldots, b_n)}(x_1,...,x_n)
  \\
  =
  \xi_{i_1}^{b_{i_1}a_{i_1}} \cdot
  \ldots \cdot
  \xi_{i_t}^{b_{i_t}a_{i_t}} \cdot
  \sum_{x_{l_1}=0}^{q_{l_1}} \xi_{l_1}^{b_{l_1}x_{l_1}} \cdot
  \ldots \cdot
  \sum_{x_{l_{s}}=0}^{q_{l_{s}}} \xi_{l_{s}}^{b_{l_{s}}x_{l_{s}}}.
 \end{multline}
 Since the number of nonzeros
 is larger than~$t$, for at least one~$j$
 from $\{l_1,\ldots,l_{s}\}$ we have
 $b_j\ne 0$.
 It follows that at least one sum in~\eqref{eq:sums} equals~$0$,
 which proves~(*).
 \end{itemize}

From (*) we conclude that
the sum of $f_C$
over $V_{i_1 ,\ldots, i_t}^{a_1 ,\ldots, a_t}$
equals the sum of~$\phi_0$ over $V_{i_1 ,\ldots, i_t}^{a_1 ,\ldots, a_t}$.
Since $\phi_0$ is a constant, the sum
does not depend on the choice
of $a_1$, \ldots,~$a_t$.

\emph{Only if.} We need to show
that $f_C$ is orthogonal
to the $\theta_i$-eigenspace for
every~$i$ in~$\{1,\ldots,t\}$,
i.e., to all
$\chi_{(b_1,...,b_n)}$
with more than~$0$ and at most~$t$
nonzeros in $(b_1,\ldots,b_n)$.
Let all $i$ such that $b_i\ne 0$
lie in $\{i_1,\ldots,i_t\}$,
where $1\le i_1<\ldots<i_t\le n$.
For the standard scalar product
$( f_C,   \chi_{(b_1,...,b_n)})$ of
$ f_C$ and $\chi_{(b_1,...,b_n)}$,
we have
\begin{multline*}
 ( f_C,   \chi_{(b_1,...,b_n)}) =
 \sum_{(x_1,...,x_n)\in V} \chi_{(b_1,...,b_n)}(x_1,...,x_n)\cdot f_C(x_1,...,x_n)
 \\
 \stackrel{\text{(i)}}=
 \sum_{(x_1,...,x_n)\in V}
 \xi_{i_1}^{b_{i_1}x_{i_1}}
 \cdot  \ldots \cdot
 \xi_{i_t}^{b_{i_t}x_{i_t}}
 \cdot f_C(x_1,...,x_n)
\\
 = \sum_{a_1=0}^{q_{i_1}-1}
 \ldots
 \sum_{a_t=0}^{q_{i_t}-1}
 \sum_{(x_1,...,x_n)\in  V_{i_1 ,\ldots, i_t}^{a_1 ,\ldots, a_t}}
 \xi_{i_1}^{b_{i_1}x_{i_1}}
 \cdot  \ldots \cdot
 \xi_{i_t}^{b_{i_t}x_{i_t}}
 \cdot f_C(x_1,...,x_n)\\
 = \sum_{a_1=0}^{q_{i_1}-1} \xi_{i_1}^{b_{i_1}a_{1}}
 \cdot
 \ldots \cdot
 \sum_{a_t=0}^{q_{i_t}-1}
 \xi_{i_t}^{b_{i_t}a_{t}}
 \cdot
 \sum_{(x_1,...,x_n)\in  V_{i_1 ,\ldots, i_t}^{a_1 ,\ldots, a_t}}
  f_C(x_1,...,x_n)\\
 \stackrel{\text{(ii)}}= \sum_{a_1=0}^{q_{i_1}-1} \xi_{i_1}^{b_{i_1}a_{1}}
 \cdot
 \ldots \cdot
 \sum_{a_t=0}^{q_{i_t}-1}
 \xi_{i_t}^{b_{i_t}a_{t}}
 \cdot
 |C|/q_{i_1}...q_{i_t}
 \stackrel{\text{(iii)}}= 0,
\end{multline*}
where in equality~(i) we use that
$\xi_j^{b_jx_j}=1$ if $b_j=0$,
in~(ii) we use the definition of an orthogonal array,
and (iii) holds because
for some $j$
in $\{i_1,...,i_t\}$ we have $b_j\ne 0$ and hence $\sum_{a=0}^{q_j}\xi_j^{b_ja}=0$.
We have shown that $f_C$ is orthogonal to $\chi_{(b_0,...,b_n)}$ for any
$(b_0,...,b_n)$ with more than~$0$
and less than $t+1$ nonzeros.
Hence, $C$ is an algebraic $t$-design.
\end{proof}

Now, we can apply the following lower bound
on the size of an algebraic design in a regular graph.

\begin{lemma}[{\cite[Sect.\,4.3.1]{KroPot:CRC&EP}}]
\label{l:gBF}
The cardinality of an algebraic $t$-design $C$
in a $k$-regular graph
$G=(V,E)$ with eigenvalues
$k=\theta_0>\theta_1>\ldots>\theta_d$ satisfies the inequality
\begin{equation}\label{eq:gBF}
 \frac{|C|}{|V|}
 \ge
 \frac{-\theta_{t+1}}{k-\theta_{t+1}}.
\end{equation}
Moreover, a multiset~$C$
of vertices of~$G$ is an algebraic $t$-design meeting~\eqref{eq:gBF} with equality if and only if $C$ is a simple set (without multiplicities more than~$1$) and a $\{k;-\theta_{t+1}\}$-CR code.
\end{lemma}
The following bound was proved
in~\cite{Friedman:92} for
$q_1=\ldots=q_s=2$ and
in~\cite{Bierbrauer:95}
 for
$q_1=\ldots=q_s=q$ for any~$q$. The theorem generalizes the mentioned results to the case of mixed-level
orthogonal arrays.
\begin{theorem}\label{th:gBF}
For an orthogonal array $C$ with parameters OA$(N,q_1{\cdot}...{\cdot}q_n,t)$,
we have
\begin{equation}\label{eq:qBF}
N \ge q_1q_2...q_n\Big(1 - \big(1-\frac1q\big)\frac{n}{t+1}\Big),
\end{equation}
where $q$
is the harmonic mean of all
$q_i$, i.e.,
$\displaystyle\frac 1q =\frac 1{n}
\sum_{i=1}^n\frac{1}{q_i}$.

Moreover, a multiset of vertices
of the graph
$H= H(q_1\cdot \ldots\cdot q_n)
$ is an OA$(|C|,q_1{\cdot}...{\cdot}q_n,t)$ meeting~\eqref{eq:qBF} with equality if and only if $C$ is a simple set (without multiplicities more than~$1$) and a $\{k;-\theta_{t+1}\}$-CR code,
where $k={Q}n(1-\frac{1}{q})$ (the degree of~$H$) and $\theta_{t+1}=k-(t+1)Q$,
$Q=\lcm(q_1, \ldots, q_n)$.
\end{theorem}
\begin{proof}
Taking into account Theorem~\ref{th:OAtD}
and Lemma~\ref{l:gBF},
it remains to check that
\eqref{eq:gBF} and~\eqref{eq:qBF} are the same for $G=H$.
Indeed, the degree of $H$ is
$$
k
=
\frac{Q}{q_1}(q_1-1) + \ldots + \frac{Q}{q_n}(q_n-1)
= Qn - Q\Big(\frac{1}{q_1}+\ldots+\frac{1}{q_n}\Big) = Qn \Big(1-\frac{1}{q}\Big)
$$
and $\theta_i=k - iQ$.
So, we find
\begin{equation*}
\frac{-\theta_{t+1}}
{k-\theta_{t+1}}=
\frac{(t+1)Q-k}
{(t+1)Q}
= 1 - \Big(1-\frac 1q\Big)\frac{n}{t+1}.
\end{equation*}
\end{proof}

Next, we compare the new bound with
the previous generalization \eqref{eq:BFD}
of~\eqref{eq:BFpu} to mixed-level orthogonal arrays~\cite{Diestelkamp:2004}.
\begin{proposition}
\label{p:better}
For mixed-level orthogonal arrays,
bound~\eqref{eq:qBF} is tighter
than~\eqref{eq:BFD}.
\end{proposition}
\begin{proof}
Rewriting~\eqref{eq:BFD} in a convenient form, we have to prove that
$$
q_1q_2...q_n\Big(1 - \big(1-\frac1q\big)\frac{n}{t+1}\Big)
\ge
q_{\Min}^n \Big(1 - \big(1-\frac1{\tilde q}\big)\frac{1}{1-(1-\frac{t+1}{n})\frac{q_{\Max}}{\tilde q}}\Big),
$$
where
\begin{align}
q_{\Min}\le q \le \tilde q \le q_{\Max}, \label{eq:mhaM}\\[1ex]
1-\frac{t+1}{n} \ge 0,\label{eq:t1n}\\[1ex]
1-\left(1-\frac{t+1}{n}\right)\frac{q_{\Max}}{\tilde q} \ge 0\label{eq:pos}
\end{align}
(\eqref{eq:mhaM} are known inequalities
between the minimum, the harmonic mean, the arithmetic mean, and the maximum values; \eqref{eq:t1n} is from $t<n$ for nontrivial orthogonal arrays;
\eqref{eq:pos} means that the denominator in~\eqref{eq:BFD}
is positive, which is the condition
of the applicability of~\eqref{eq:BFD}).
The required inequality is straightforward from the following three observations.
\begin{itemize}
    \item[(i)] Trivially, $q_1q_2...q_n \ge q_{\Min}^n$.

\item[(ii)] From~\eqref{eq:mhaM} we have $\displaystyle 1-\frac1q  \le 1-\frac1{\tilde q}$.

\item[(iii)] Taking into account \eqref{eq:mhaM}--\eqref{eq:pos}, we get
$$
\frac{n}{t+1}
=
\frac{1}{1-(1-\frac{t+1}{n})}
\le
\frac{1}{1-(1-\frac{t+1}{n})\frac{q_{\Max}}{\tilde q}}.
$$
\end{itemize}

Finally, each of (i), (ii), (iii)
turns to equality if and only if
$q_1=q_2=\ldots=q_n$, i.e., for pure-level arrays.
\end{proof}

\begin{remark}
As noted in~\cite{Diestelkamp:2004},
the definition of an orthogonal
array implies that
$N$ must be a multiple
of $$Q_t = \lcm\,\{q_{i_1}q_{i_2}\ldots q_{i_t}:\ 0 < i_1 < \ldots <i_t\le n  \}.$$
In particular, this gives the bound
$N\ge Q_t$ (by similarity of arguments, it can be considered as an analog of the Singleton bound for
error-correcting codes),
but also means that any bound
of form $N\ge B(q_1\cdot\ldots\cdot q_n,t)$ can be rounded to
\begin{equation}\label{eq:round}
 N\ge \lceil B(q_1\cdot\ldots\cdot q_n,t)/Q_t  \rceil \cdot Q_t.
\end{equation}
However, in some cases bound \eqref{eq:qBF} is already divisible
by~$Q_t$, and the problem of the existence
of orthogonal arrays attaining it arises.
In the next section, we show how to construct such arrays in the case when all $q_i$ are powers of the same prime number.
\end{remark}

%%%%%%%%%%%%%%%%%%%%%%%%%%%
\section{Additive mixed-level orthogonal arrays attaining the Bierbrauer--Friedman bound}\label{s:additive}

In this section, all $q_i$
are powers of the same prime~$p$,
and $\intrv{p^i}$ is associated
with the vector space $\FF_p^i$
(to be explicit, an integer $a$ in $\intrv{p^i}$ can be associated
with its $p$-based notation,
treated as a vector from $\FF_p^i$).
We should warn the reader about the following difference in notation:
now $q_i$ denotes $p^i$, not the number of levels
(alphabet size) in the $i$th position.
This is because in this section,
it is convenient to write parameters
in the form OA$(N,q_1^{n_1}q_2^{n_2}\ldots q_s^{n_s},t)$
(all $q_i$ are different, but some $n_i$, $i<s$, can be zero), while in Section~\ref{s:OA-AD}
the preferred form was
OA$(N,q_1\cdot q_2\cdot \ldots\cdot q_n,t)$, where $q_i$ are not necessarily distinct.

\begin{definition}[additive code]
A set of $(n_1+\ldots+n_s)$-tuples from
$V=(\FF_p^1)^{n_1}\times (\FF_p^2)^{n_2} \times \ldots \times (\FF_p^s)^{n_s} $,
where $p$ is prime,
is called \emph{additive}
(an additive code, or an additive orthogonal array if we consider it as an orthogonal array)
if it is closed with respect to
the coordinatewise addition, i.e.,
forms a subspace of the $(n_1+2n_2+\ldots+sn_s)$-dimensional vector space~$V$
over~$\FF_p$.
\end{definition}

\begin{remark}
For prime $p$, ``additive'' and ``$\FF_p$-linear'' are the same, but
in general, if $p$ is a prime power,
$\FF_p$-linear codes form a proper
subclass of additive codes.
The theory in the rest
of this section keeps
working for an arbitrary prime power~$p$
if we replace ``additive'' by ``$\FF_p$-linear'' everywhere. However, to simplify reading, we focus on the most
important case of prime~$p$ and localize the arguments for the general case in this remark.
\end{remark}

Any $k$-dimensional vector subspace~$C$ of an $n$-dimensional
vector space~$S_n$ can be represented as the null-space (kernel)
of a homomorphism from~$S_n$
to an $(n-k)$-dimensional vector space~$S_{n-k}$ over the same field.
If the bases in $S_n$ and $S_{n-k}$, are fixed,
such a homomorphism is represented
by an $(n-k)\times n$ matrix,
called a \emph{check matrix} of~$C$.
In our case the dimension of the space is $(n_1+2n_2+\ldots+sn_s)$,
we have a natural basis,
and the coordinates of a vector,
as well as the columns of a check matrix,
are naturally divided into $n_1+n_2+\ldots+n_s$ groups,
\emph{blocks}, the first $n_1$ blocks of size~$1$, the next $n_2$ blocks of size~$2$, and so on.
For a check matrix $H$,
by $H_{i,j}$ we denote the space
spanned by the columns from the
$j$th block of size~$i$.

\begin{definition}[multifold partition of a space]\label{d:mu}
A multiset $D$ of subspaces of a vector space~$S$
is called a \emph{$\mu$-fold partition} of~$S$ if
every nonzero vector of~$S$ belongs to exactly $\mu$
subspaces from~$D$, respecting the multiplicities.
\end{definition}
Next, we define two concepts,
which are not (in contrast
to Definition~\ref{d:mu}) key
concepts in our theory, but allow
to mention additionally one important correspondence in the main theorem of this section.
\begin{definition}[one-weight code, alphabet-effective code]
An additive code $C$ in $V$ is called
\emph{one-weight} (of weight~$w$)
if the number of nonzero blocks
in a nonzero codeword is constant
(equal to~$w$).
The code is called \emph{alphabet-effective}
if in each block, every element of the corresponding alphabet $\FF_p^i$
occurs in some codeword
(so, the alphabet $\FF_p^i$ is effectively used).
\end{definition}

\begin{theorem}\label{th:additive}
Let $H$ be an $m\times (n_1+2n_2+\ldots+sn_s)$ rank-$m$ matrix over~$\FFp$, and let $C$ be the null-space of~$H$ (so, $C$ is an additive code in $(\FF_p^1)^{n_1}\times \ldots \times (\FF_p^s)^{n_s} $). Denote
\begin{equation}\label{eq:dnt}
  k=\sum_{i=1}^s n_i(p^s-p^{s-i})
  =p^s\sum_{i=1}^s n_i(1-p^{-i})
  ,
\quad
\mu = \frac{k}{p^m-1},
\quad q_i=p^i.
\end{equation}
The following assertions
are equivalent:
\begin{itemize}
\item[\emph{(i)}] $C$ is an
OA$(|C|,q_1^{n_1}q_2^{n_2}...q_s^{n_s},t)$ attaining bound~\eqref{eq:qBF};
\item[\emph{(i')}] $C$ is an
OA$(|C|,q_1^{n_1}q_2^{n_2}...q_s^{n_s},t)$
with $t= \frac{\mu+k}{p^s}-1$;
\item[\emph{(ii)}] $C$ is a $\{k,\mu\}$-CR code in $H(q_1^{n_1}q_2^{n_2}...q_s^{n_s})$;
\item[\emph{(iii)}] the code $C^\perp$ generated by the rows of $H$ is an alphabet-effective one-weight code of weight~$\mu p^{m-s}$;
\item[\emph{(iv)}] the multiset
$M=\{p^{s-i} \times H_{i,j}\}_{i=1}^s\vphantom{)}_{j=1}^{n_i} $
of subspaces of~$\FFp^m$ is a  $\mu$-fold
partition ($p^{s-i} \times H_{i,j}$ denotes that $H_{i,j}$ is added $p^{s-i}$ times in the multiset);
\item[\emph{(v)}] the collection
$M^{\perp}=\{H_{i,j}^{\perp}\}_{i=1}^s\vphantom{)}_{j=1}^{n_i} $
of subspaces dual to~$H_{i,j}$
is a $\nu$-fold
partition of~$\FFp^m$,
where $\nu = n_1+\ldots+n_s - p^{m-s}\cdot \mu$.
%for every $i\in\{1,...,s\}$ and $j\in\{1,...,n_i\}$ the subspace $H_{i,j}$ has rank~$i$ and
%every nonzero vector
%of~$\FFp^m$ is exactly
%$-\theta_{t+1}=-\sum_{i=1}^s n_i(p^s-p^{s-i})+p^s(t+1)$ times covered by subspaces from the multiset $\displaystyle \biguplus_{i=1}^s \biguplus_{j=1}^{n_i} \{\underbrace{H_{i,j},\ldots,H_{i,j}}_{p^{s-i}\text{ times}}\}$.
\end{itemize}
\end{theorem}
\begin{proof}
Since $C$ is a null-space of $H$,
we have $|C| = p^{n_1+2n_2+\ldots+sn_s}/p^m =
q_1^{n_1}\cdot
q_2^{n_2}\cdot \ldots \cdot q_s^{n_s} / p^m$.

(i)$\Longleftrightarrow$(i')
Taking into account the expression for $|C|$ above, \eqref{eq:qBF} has the form
\begin{equation}\label{eq:BFp}
p^{-m} \ge \Big(1-\big(1-\frac1q\big)\frac{n}{t+1}\Big).
\end{equation}
(recall that $q$ is the harmonic mean of
$
\underbrace{q_1,...,q_1}_{n_1\text{ times}},
\underbrace{q_2,...,q_2}_{n_2\text{ times}},
\ldots,
\underbrace{q_s,...,q_s}_{n_s\text{ times}}
$).
From the equality in~\eqref{eq:BFp},
we find
$$
t+1=\frac{n(1-\frac1q)}{1-p^{-m}}
= \frac{p^{-s} k}{1-p^{-m}} =
p^{-s} k \Big(1+\frac{1}{p^m-1}\Big) = \frac{k+\mu}{p^s}.
$$
Inversely, substituting $t = \frac{k+\mu}{p^s}-1$ turns \eqref{eq:BFp}
to equality.

(i)$\Longleftrightarrow$(ii)
By
Theorem~\ref{th:gBF},
$C$ attains~\eqref{eq:qBF}
if and only if it is a
$\{k;-\theta_{t+1}\}$-CR code,
where $k$ is the degree of
$H(q_1^{n_1}q_2^{n_2}...q_s^{n_s})$
and $\theta_{t+1} = k-p^s(t+1)$
(the $(t+1)$th largest eigenvalue of $H(q_1^{n_1}q_2^{n_2}...q_s^{n_s})$, counting from~$0$).

(ii)$\Longleftrightarrow$(iv)
Assume that (iv) holds
and we have to show~(ii).
We first check that

\begin{itemize}
    \item[(*)] \emph{$H_{i,j}$
has dimension~$i$ for each~$i$, $j$.}
\end{itemize}

Indeed, if it is so, then $H_{i,j}$
has $(p^i-1)$ nonzero points,
$p^{s-i}\times H_{i,j}$ has
$(p^s-p^{s-i})$ nonzero points,
and $M$ has
$\sum_{i=1}^s n_i(p^s-p^{s-i})$,
i.e., $k$, nonzero points.
This is exactly the number of points
we need to cover all $p^m-1$
nonzero points of~$\FF_p^m$
with multiplicity
$\mu = \frac{k}{p^m-1}$.

But if one $H_{i,j}$ has dimension smaller than $i$, then the number
of points in~$M$ is not enough
to make a $\mu$-fold partition
of~$\FF_p^m$. So, (*) holds.

Next, (*) means that $C$ is independent.
Indeed, if there are two
adjacent
vectors $x$ and~$y$ in~$C$,
then their difference $x-y$
has nonzero values only in coordinates from one block, say $(i,j)$th.
In this case, these values are
coefficients of a nontrivial
linear dependency between the corresponding columns,
and hence $H_{i,j}$ has dimension smaller than~$i$, contradicting~(*).

The independence of $C$ implies
that every its element is connected
to elements not in~$C$ by $k$ edges
($k$ is the degree of the graph),
and it remains to confirm the second
parameter~$\mu$ of the $\{k;\mu\}$-CR code.

For a vector $v$ not in $C$
denote by $s$ its \emph{syndrome}
$H\cdot v^{\mathrm T}$ (since $c\not\in C$, the syndrome $s$ is nonzero).
A vector $v-e$ is adjacent to~$v$ and belongs to~$C$
if and only if $H\cdot e^{\mathrm T}=s$
and
$e$ has zeros in all positions out of
one block, say~$(i,j)$th.
The last is equivalent to $s\in H_{i,j}$.
Since the subspace $H_{i,j}$ is counted
with multiplicity $p^{s-i}$ in the multiset~$M$, and the multiplicity
of the edge $\{v,v-e\}$ is also $p^{s-i}$,
the number of edges between $v$ and~$C$ is equal to the multiplicity of~$s$ in~$M$,
i.e.,~$\mu$.

By reversing the arguments, we ensure that
(ii) implies~(iv).

(iv)$\Longleftrightarrow$(v)
The equivalence between (iv) and (v)
is proven in~\cite[Theorem~3]{KroMog:multispread}
(which is, apart of the direct correspondence between the parameters,
is essentially a special case of
\cite[Theorem~15]{El-Zanati-et-al:lambda}).

(iv)$\Longleftrightarrow$(iii)
The equivalence between (iv) and (iii)
is essentially proven in~\cite[Theorem~1]{KroMog:multispread}.
The difference is that \cite[Theorem~1]{KroMog:multispread} considers one-weight codes over the same alphabet $\FF_p^s$ (in our current notation), but not necessarily alphabet-effective.
It remains to observe that there is a trivial weight-preserving relation between alphabet-effective additive codes over our mixed alphabet and additive codes over the alphabet $\FF_p^s$, namely, adding $s-i$ zero columns to each block of size $m\times i$
of the generator matrix~$H$.
\end{proof}

A collection of subspaces of $\FFp^m$
of dimension at most~$s$
satisfying the condition of Theorem~\ref{th:additive}
is called a $(\lambda,\mu)^{s,m}_p$-multispread \cite{KroMog:multispread}, where
% $\mu=-\theta_{t+1}$ and
$\lambda=n_1 (p^{s-1}-1)+n_2 (p^{s-2}-1)+\ldots+n_{s-1} (p^{1}-1)$
(the multiplicity of the zero vector in~$M$).
Theorem~\ref{th:additive}
means that such multispreads with $n_i$ subspaces of dimension $i$, $i=1,...,s$,
corresponds to orthogonal
arrays in $(\FFp^1)^{n_1}\times \ldots\times (\FFp^s)^{n_s}$ attaining the generalized Bierbrauer--Friedman bound~\eqref{eq:qBF} in Theorem~\ref{th:gBF}.
Note that $n_1$, \ldots, $n_s$ are not
uniquely determined from the parameters
$p$, $s$, $m$, $\lambda$, $\mu$ of a $(\lambda,\mu)^{s,m}_p$-multispread,
and even without fixing $n_1$, \ldots, $n_s$
the problem of existence of $(\lambda,\mu)^{s,m}_p$-multispread is open in general (in~\cite{KroMog:multispread},
it is completely solved for $s=2$, any~$p$,
and for $p^s\in\{2^3,2^4,3^3\}$).
The pure-level subcase
of Theorem~\ref{th:additive}
was proved in~\cite[Theorem~4.8]{BKMTV};
in that case, $\lambda=0$, and multispreads are $\mu$-fold spreads, whose parameters are characterized, see
\cite[p.83]{Hirschfeld79}, \cite[Corollary~8]{El-Zanati-et-al:lambda}.

\begin{example}\label{ex:1}
Consider the $3\times (1\cdot 1+4\cdot 2)$ matrix
$$H=\left(
\begin{array}{cc@{\ }cc@{\ }cc@{\ }cc@{\ }c}
    1 &  1 & 0  &  1 & 0  &  0 & 0  &  1 & 0  \\
    1 &  0 & 1  &  0 & 0  &  1 & 0  &  1 & 1  \\
    1 &  0 & 0  &  0 & 1  &  0 & 1  &  0 & 1
\end{array}\right )$$
over~$\FF_2$. Its null-space $C$ is an additive
OA$(64,2^14^4,3)$ attaining~\eqref{eq:qBF}.
Indeed, 
\begin{itemize}
    \item $64 = 2^9\big(1-(1-\frac{3}{10})\frac{5}{3+1}\big)$;
    \item the degree of the graph $H(2^1 4^4)=2K_2\times K_4\times K_4\times K_4\times K_4$ is $2+4\cdot 3=14$;
    \item the eigenvalue $\theta_{t+1}=\theta_{4}$ equals
$14 - 4\cdot 4=-2$;
    \item the subspace
$H_{1,1}=\langle (1,1,1) \rangle$,
corresponding to the first column, is taken with multiplicity~$2$ and covers the vector $(1,1,1)$ exactly $2=-\theta_{4}$ times; each of the other nonzero vectors
from~$\FF_2^3$ is covered exactly $2=-\theta_{4}$ times by the subspaces
$H_{2,1}=\langle (1,0,0),(0,1,0) \rangle$,
$H_{2,2}=\langle (1,0,0),(0,0,1) \rangle$,
$H_{2,3}=\langle (0,1,0),(0,0,1) \rangle$,
$H_{2,4}=\langle (1,1,0),(0,1,1) \rangle$, formed by the last four blocks of columns;
so, $\mu=2$.
    \item According to p.(iv) of Theorem~\ref{th:additive},
the dual subspaces
$H_{1,1}^\perp$, $H_{2,1}^\perp$, $H_{2,2}^\perp$, $H_{2,3}^\perp$, $H_{2,4}^\perp$
form a partition of~$\FF_2^3$,
$\nu=n_1+n_2-\mu\cdot 2^{m-s}=1+4-2\cdot 2^{3-2}=1$.
\item Every nonzero linear combination 
of the rows of~$H$ has exactly $4$ nonzero blocks, i.e, $C^\perp$ is a one-weight code of weight $4 = \mu\cdot 2^{m-s}$.
\end{itemize}

\end{example}

% \alert{characterize for $p^s=p^2$, other?}

\section{Polynomial generalization of the Bierbrauer--Friedman bound }\label{s:PG-BFB}

For a connected graph $G$, $\diam(G)$
denotes its diameter, and $G^{(j)}$,
denotes a simple graph where two vertices are
adjacent if and only if the distance between them in~$G$ is~$j$.

\begin{definition}[{distance-regular graph}]
\label{1def:drg}
A simple connected graph
is called a \emph{distance-regular graph}
 if for $i = 1,\ldots, \diam(G)$ the product $A^{(1)}A^{(i)}$ is a linear combination of $A^{(i-1)}$, $A^{(i)}$, and $A^{(i+1)}$,
 where $A^{(j)}$ is the adjacency matrix of~$G^{(j)}$.
\end{definition}

\begin{corollary} \label{1c:drgP}
 For every distance-regular graph~$G$,
 there are degree-$i$ polynomials $K^{(i)}$, $i=0,\ldots,\diam(G)$,
 such that $G^{(i)} = K^{(i)} (G)$.
\end{corollary}

Let $G$ be a distance-regular graph.
Consider a nonnegative integer linear combination of the polynomials $K^{(i)}$, i.e., $P=\sum_i \alpha_iK^{(i)}$,
$a_i\geq 0$. Then $P(A)$ is the adjacency matrix of a multigraph $P(G)$ with the same set~$V$ of vertices. If vertices~$u$ and~$v$ are at distance~$i$ from each other in~$G$, then the edge $\{u,v\}$ has multiplicity~$\alpha_i$ in~$P(G)$. The following proposition collects some obvious and straightforward facts about eigenfunctions of $G$ and~$P(G)$.

\begin{proposition}\label{p:gBF25}
{\emph{(a)}} The set of eigenfunctions of $P(G)$ includes  the set of eigenfunctions of~$G$.
If
$P$
is strictly monotonic function, then these sets coincide.
\\
{\emph {(b)}} If $\varphi$ is an eigenfunction of $G$ with eigenvalue~$\theta$,
then $\varphi$ is also an eigenfunction of
$P(G)$
with eigenvalue
$P(\theta)$.\\
{\emph {(c)}} If  $P$  strictly increases, then the sets of algebraic $t$-designs of~$G$ and~$P(G)$ coincide.
\end{proposition}

The next lemma is a variant of Lemma~\ref{l:gBF}.

\begin{lemma}
\label{l:gBF1}
Let $G=(V,E)$ be a distance-regular graph with eigenvalues
$k=\theta_0,\theta_1,\ldots,\theta_{\diam(G)}$, in the decreasing order, let the polynomial~$P$ be a linear combination
$$ P(x) = \alpha_1 K^{(1)}(x)+\ldots +
\alpha_{\diam(G)} K^{(\diam(G))}(x)
$$
of the polynomials~$K^{(i)}$
 from Corollary~\ref{1c:drgP},
 where $\alpha_i\ge 0$,
 $i=1,\ldots,\diam(G)$,
 and let $$\mu = \max\, \{ P(\theta_{t+1}),\ldots,P(\theta_{\diam(G)})\}.$$
If $\mu<0$, then the cardinality of an algebraic $t$-design $C$ satisfies the inequality
\begin{equation}\label{eq:pBF1}
 \frac{|C|}{|V|}
 \ge
 \frac{-\mu}{P(k)-\mu}.
\end{equation}
\end{lemma}
\begin{proof}
If all $\alpha_i$ are integer,
then $P(x)$ is a graph,
for which Lemma~\ref{l:gBF} gives
\eqref{eq:pBF1}.

If all $\alpha_i$ are rational,
we divide $P$ by their greatest common divisor, which makes all $\alpha_i$ integer and does not affect the value in~\eqref{eq:pBF1}.

If some $\alpha_i$ are irrational,
we can approximate them by rational values
with any accuracy
and get~\eqref{eq:pBF1} as the limit
of bounds corresponding to polynomials~$P$ with rational coefficients.
\end{proof}

\begin{remark}
    Although
in this section we focus on distance-regular graphs, the approach is partially applicable in more general cases. Proposition~\ref{p:gBF25} and Lemma~\ref{l:gBF1} above work for any regular graph~$G$ and a polynomial~$P$ such that for the adjacency matrix~$A$ of~$G$, the matrix~$P(A)$ is nonnegative (and integer, for Proposition~\ref{p:gBF25}). Combining the distance-$i$ polynomials $K^{(i)}$ of~$G$ is also possible if they exist for some values of~$i$, not necessarily for all~$i$ from $0$ to~$\mathrm{d}(G)$. \\ \end{remark}
Our aim is to find a polynomial~$P$ such that $\theta_{t+1}\geq 0$ but $\mu<0$.
In this case, the lower bound from
Lemma~\ref{l:gBF} is non-positive and hence trivial,
while the lower bound from Lemma~\ref{l:gBF1} is positive and can be valuable.
Actually, to obtain the best lower bound in~\eqref{eq:pBF1},
we need to solve the following linear programming problem in the variables $\alpha_1$, \ldots,
$\alpha_{\diam(G)}$, $\mu$
(the last one is an additional variable, which plays the role of $\max \{P(\theta_{j})\}_{j\ge t+1}$):
$$
 \left\{
\begin{array}{l}
 \sum\limits_{i=1}^{\diam(G)}\alpha_i
 K^{(i)}(\theta_j)\le \mu, \qquad  j = t+1,\ldots,\diam(G),\\[3ex]
\alpha_i\geq 0,\qquad i= 1,\ldots,\diam(G),\\[3ex]
 \sum\limits_{i=1}^{\diam(G)}a_i K^{(i)}(\theta_{0})=\mathrm{const}
 \quad \mbox{(the degree of $P(G)$)},\\[3ex]
\mu \longrightarrow \min.
\end{array}
\right.
$$

By the proof of Lemma~\ref{l:gBF1}, we do not need to care about the integrity of $\alpha_i$.
By Lemmas \ref{l:gBF}, \ref{l:gBF1}, and Proposition \ref{p:gBF25}(c) we can conclude the following.
\begin{proposition}
Let a collection $a_1,\dots,a_{\diam(G)}$ be the solution of the  linear program and let $P$ be the corresponding polynomial. Suppose that $P(\theta_{t+1})\neq P(\theta_{j})$ for $j\neq t+1$.  If some $t$-design~$C$ attains the bound \eqref{eq:pBF1}, then $C$ is a CR-$1$ code in~$G$.
\end{proposition}

For the Hamming graph $G=H(n,q)$, the degree-$w$ polynomial~$K^{(w)}$
 is
  \begin{multline*}\label{1eq:krawtchuk}
  K^{(w)}(\cdot) = P_w(P_1^{-1}(\cdot)),
  \\ \mbox{where }
   P_w(x) = P_w(x;n,q) = \sum_{j=0}^w (-1)^j (q-1)^{w-j}
   \Big(\genfrac{}{}{0mm}{}{x}{j} \Big)
   \Big(\genfrac{}{}{0mm}{}{n-x}{w-j} \Big)
  \end{multline*}
 is the \emph{Krawtchouk polynomial}; $P_1(x) = (q-1)n-qx$, \ $\displaystyle P_1^{-1}(y)=\frac{(q-1)n-y}q $.

\begin{example}
Consider $G=H(2m,2)$, $t=m-1$.
Let us derive a bound on the size
of OA$(N,2m,2,m-1)$
by combining only three polynomials
$K^{(i)}$, $i=1,2,3$.
We have $G=H(2m,2)$, $\theta_{t+1}=\theta_{m}=0$,
$$K^{(1)}(x)=x,\quad
K^{(2)}(x)=\frac12(x^2-2m),
\quad
K^{(3)}(x)=\frac16\big(x^3-(6m-2)x\big).$$
We are looking for a polynomial $P=K^{(3)}+\beta K^{(2)}+\alpha K^{(1)}$ where $\alpha\ge 0$, $\beta\ge 0$, and $P'(x)\geq 0$ for all $x\in [-2m,0]$.
Put $\alpha=\frac{\beta^2}{2}+m-\frac{1}{3}$. Then $6P(x)=(x+\beta)^3-\beta^3-6\beta m$. It is easy to see that
$P$ increases monotonically. Substituting  $\beta=\frac{m}{\sqrt{3}}$, we obtain
$6P(x)
= x^3 + \sqrt3 mx^2 + m^2 x -2\sqrt3 m^2$ and, by~\eqref{eq:pBF1},
\begin{multline*}
\frac{|C|}{|V|}\geq \frac{-6P(0)}{6P(2m)-6P(0)}
=\frac{2\sqrt 3 m^2}{10m^3 + 4\sqrt3 m^3 }
=
\frac{1}{m(5\cdot 3^{-\frac12}+2)}>\frac{0.2046}{m}.
\end{multline*}
\end{example}

\section{Conclusion}
The most general contribution
of our correspondence
is a representation of
orthogonal arrays as algebraic designs
(which was previously known only for pure-level arrays).
More specifically, as an application of
this representation,
we have considered two possibilities to generalize the Bierbrauer--Friedman bound for orthogonal arrays.

One of the generalizations
concerns mixed-level
orthogonal arrays,
and we have also shown how to construct
arrays attaining the generalized bound.
However, the constructed arrays have
one restriction because of their additive structure over a finite field;
namely, the size of such an array
is always a prime power.
As was mentioned in the introduction,
the only known arrays attaining the Bierbrauer--Friedman bound with the size
not being a prime power
are pure-level binary or ternary. By generalizing the bound the problem of the existence of such arrays also expands to the mixed-level case. In particular: do there exist mixed-level OA$(N,2^{n_1}3^{n_2},t)$ attaining bound~\eqref{eq:qBF}? Orthogonal arrays attaining the considered bound (in both pure- and mixed-level cases) are of special interest because of their relation with such regular structures as CR-$1$ codes (equivalently, intriguing sets, equitable $2$-partitions, perfect $2$-colorings).

The other, polynomial, generalization,
considered in Section~\ref{s:PG-BFB},
makes it possible to obtain a lower bound
on the size of a pure-level orthogonal array when the original Bierbrauer--Friedman bound gives a trivial inequality.
The authors do not know if it is possible
to improve the Bierbrauer--Friedman bound in this way in cases when it is positive.

\section*{Declaration of competing interest, data availability}
\ldots no known competing interests \ldots
no data were used or generated \ldots.

%%%%%%%%%%%%%%%%%%%%%%%%%%%%%%%%%%%%%%%%
% \bibliographystyle{plain}
% \bibliography{k}
% \end{document}

\providecommand\href[2]{#2} \providecommand\url[1]{\href{#1}{#1}}
  \def\DOI#1{{\href{https://doi.org/#1}{https://doi.org/#1}}}\def\DOIURL#1#2{{\href{https://doi.org/#2}{https://doi.org/#1}}}

\end{document}